\documentclass[10pt]{article}%
\usepackage{makeidx}
\usepackage{graphicx}
\usepackage{amsmath}
\usepackage{amsfonts}
\usepackage{amssymb}%
\providecommand{\U}[1]{\protect\rule{.1in}{.1in}}
\newtheorem{theorem}{Theorem}

\newtheorem{corollary}[theorem]{Corollary}

\newtheorem{definition}[theorem]{Definition}
\newtheorem{example}[theorem]{Example}

\newtheorem{proposition}[theorem]{Proposition}

\newenvironment{proof}[1][Proof]{\textbf{#1.} }{\ \rule{0.5em}{0.5em}}
\begin{document}

\title{Polyharmonicity and algebraic support of measures}
\author{Ognyan Kounchev and Hermann Render}
\maketitle

\begin{abstract}
Our main result states that two measures $\mu$ and $\nu$ with bounded support
contained in the zero set of a polynomial $P(x)$ are equal if they coincide on
the subspace of all polynomials of polyharmonic degree $N_{P}$ where the
natural number $N_{P}$ is explictly computed by the properties of the
polynomial $P\left(  x\right)  $. The method of proof depends on a definition
of a multivariate Markov transform which another major objective of the
present paper. The classical notion of orthogonal polynomial of second kind is
generalized to the multivariate setting: it is a polyharmonic function which
has similar features as in the one-dimensional case.

\textrm{Acknowledgement}: \emph{Both authors have been sponsored by the
Institutes partnership project at the Alexander von--Humboldt Foundation; the
first was sponsored by a Greek--Bulgarian S\&T\ Cooperation project.}



\end{abstract}

\section{Introduction}

Recall that a complex-valued function $f$ defined on a domain $G$ in the
euclidean space $\mathbb{R}^{n}$ is \emph{polyharmonic of order }$N$ if $f$ is
$2N$-times continuously differentiable and
\[
\Delta^{N}f\left(  x\right)  =0\text{ for all }x\in G
\]
where $\Delta^{N}$ is the $N$-th iterate of the Laplace operator $\Delta
=\frac{\partial^{2}}{\partial x_{1}^{2}}+...+\frac{\partial^{2}}{\partial
x_{n}^{2}}.$ For $N=1$ this class of functions are just the harmonic
functions, while for $N=2$ the term biharmonic function is used which is
important in elasticity theory. Fundamental work about polyharmonic functions
is due to E. Almansi \cite{Alma99}, M. Nicolesco (see e.g. \cite{Nico35}) and
N. Aronszajn \cite{ACL83}, and still this is an area of active research, see
e.g. \cite{Eden}, \cite{FKM01},\cite{FKM03}, \cite{Hayman}, \cite{kounchev92}%
,\cite{kounchev98}, \cite{Ligo88}, \cite{Rend05}, \cite{Sobo74}. Polyharmonic
functions are also important in applied mathematics, e.g. in approximation
theory, radial basis functions and wavelet analysis, see e.g. \cite{BBRV05},
\cite{Koun00}, \cite{KoReCM}, \cite{kounchevrenderJAT}, \cite{MaNe90}.

In this paper we address the following question: suppose that $P\left(
x\right)  $ is a polynomial, and that $\mu$ and $\nu$ are signed measures
which have support in the zero set $K_{P}$ of the polynomial $P,$ i.e. in the
set
\[
K_{P}\left(  R\right)  :=\left\{  x\in\mathbb{R}^{n}:P\left(  x\right)
=0\text{ and }\left|  x\right|  \leq R\right\}  .
\]
Under which conditions do $\mu$ and $\nu$ coincide? As motivating example
consider the polynomial $P\left(  x\right)  =\left|  x\right|  ^{2}-1$ where
$\left|  x\right|  :=r\left(  x\right)  :=\sqrt{x_{1}^{2}+...+x_{n}^{2}}$ is
the euclidean norm in $\mathbb{R}^{n}.$ It is well known that two measures
$\mu$ and $\nu$ with support in the unit sphere $\mathbb{S}^{n-1}=\left\{
x\in\mathbb{R}^{n}:\left|  x\right|  =1\right\}  $ coincide if they are equal
on the set of all \emph{harmonic} polynomials. We shall show that two measures
$\mu$ and $\nu$ with support in $K_{P}\left(  R\right)  $ are equal if the
moments $\mu\left(  f\right)  $ an $\nu\left(  f\right)  $ are equal for
polyharmonic polynomials $f$ of a certain degree $N_{P}$ which depends on the
polynomial $P.$ In order to formulate this precisely, let us introduce the
\emph{polyharmonic degree} $d\left(  f\right)  $ defined by
\begin{equation}
d\left(  f\right)  :=\min\left\{  N\in\mathbb{N}_{0}:\Delta^{N+1}\left(
f\right)  =0\right\}  \label{pdegree}%
\end{equation}
Note that$f$ has \textbf{polyharmonic degree} $\leq N$ if and only if $f$ is
of \textbf{polyharmonic order} $N+1.$

Let us denote by $\mathcal{P}$ set of all polynomials. One of the main results
of this paper reads as follows:

\begin{theorem}
\label{ThmMain}Let $P\left(  x\right)  $ be a polynomial and define
\[
N_{P}:=\sup\left\{  d\left(  P\cdot h\right)  :h\in\mathcal{P}\text{ is a
harmonic polynomial}\right\}  .
\]
Let $\mu$ and $\nu$ be measures with support contained in the set
$K_{P}\left(  R\right)  $ for some $R>0.$ Then $\mu\equiv\nu$ if and only if
$\int hd\mu=\int hd\nu$ for all polynomials $h$ in the subspace
\[
U_{N_{P}}:=\left\{  Q\in\mathcal{P}:\ \Delta^{N_{P}}Q=0\right\}  .
\]

\end{theorem}

It is easy to see that $N_{P}$ is lower or equal to the total degree of the
polynomial $P\left(  x\right)  .$ In the appendix we shall give a procedure to
determine the number $N_{P}$ explicitly.

An application of the Hahn-Banach theorem shows us the following consequence
of Theorem \ref{ThmMain}: the space $U_{N_{P}}$ is dense in the space
$C\left(  K_{P}\left(  R\right)  ,\mathbb{C}\right)  $ of all continuous
complex-valued functions on the compact space $K_{P}\left(  R\right)  $
endowed with the supremum norm, see Corollary \ref{Cdensity}. We call the
reader's attention to this interesting result which may be compared with the
density results for solutions to $\Delta^{p}h=0$ in $C\left(  K\right)  $ for
compacts $K,$ obtained with the techniques of Potential theory in the $1970$s;
see \cite{hedberg}, \cite{hedberg2} and the references therein.

It is also instructive to consider the statement of Theorem \ref{ThmMain} for
the univariate case $n=1$, so $P$ is a polynomial of degree $N,$ and
$P^{-1}\left(  0\right)  $ has at most $N$ elements. Note that $\Delta
^{N}Q=\frac{d^{2N}}{dx^{2N}}Q=0$ if and only if $Q$ is a polynomial of degree
$\leq2N-1.$ Hence, Theorem \ref{ThmMain} says that two non-negative measures
$\mu$ and $\nu$ with support in $P^{-1}\left(  0\right)  $ are equal if and
only if
\[
\int x^{s}d\mu=\int x^{s}d\nu\text{ for all }s\leq2N-1.
\]
So Theorem \ref{ThmMain} can be seen as a generalization of a simple
univariate statement based upon the Polyharmonic paradigm as presented in
\cite[chapter 1.5]{Koun00}.

The proof of Theorem \ref{ThmMain} will be a by-product of our investigation
of the so-called \emph{multivariate Markov transform} which we will introduce
below and which we consider as a suitable generalization of the univariate
\emph{Markov transform}, an important tool in the classical moment problem and
its applications to Spectral theory. Recall that the \emph{Markov
transform}\footnote{In some recent works in Approximation theory, Potential
theory, and Probability theory this function is called the \emph{Markov
function} of a measure, see e.g. \cite{StTo92} or \cite{gonchar}. On the other
hand apparently Widder \cite{widder} was the first who has given the name
\emph{Stieltjes transform }to this function. If $\mu$ has infinite support the
transform is also called Stieltjes transform. This tradition has been followed
by Akhiezer \cite{Akhi65} and other Russian mathematicians.} of a finite
measure $\sigma$ with support in the interval $\left[  -R,R\right]  $ is
defined on the upper half--plane by the formula
\begin{equation}
\widehat{\sigma}\left(  \zeta\right)  :=\int_{-R}^{R}\frac{1}{\zeta-x}%
d\sigma\left(  x\right)  \text{ for }\operatorname{Im}\zeta>0,
\label{defstieltjes}%
\end{equation}
see e.g. \cite[Chapter 2]{Akhi65}, \cite[Chapter 2.6]{NiSo91}. Let us recall a
central result called Markov's theorem: the $N-$th Pad\'{e} approximant
$\pi_{N}\left(  \zeta\right)  =Q_{N}\left(  \zeta\right)  /P_{N}\left(
\zeta\right)  $ of the asymptotic expansion of $\widehat{\sigma}\left(
\zeta\right)  $ at infinity converges compactly in the upper half plane to
$\widehat{\sigma}\left(  \zeta\right)  ;$ here the polynomial $P_{N}$ is the
$N$-th orthogonal polynomial with respect to the measure $\sigma$ and $Q_{N}$
is the \emph{orthogonal polynomial of the second kind} with respect to the
measure $\sigma$ given through the formula
\begin{equation}
Q_{N}\left(  \zeta\right)  =\int_{-\infty}^{\infty}\frac{P_{N}\left(
\zeta\right)  -P_{N}\left(  x\right)  }{\zeta-x}d\sigma\left(  x\right)  .
\label{secondkind}%
\end{equation}
Further, to each $\pi_{N}\left(  \zeta\right)  $ there corresponds a
(non--negative) measure $\sigma_{N}$ with support in the zeros of the
nominator $P_{N},$ thus leading to a proof of the famous Gau\ss \ quadrature formula.

Our definition of a multivariate Markov transform depends on the work of
\emph{N. Aronszajn} \cite{ACL83} on polyharmonic functions, and of \emph{L.K.
Hua} \cite{Hua63} about harmonic ana\-lysis on Lie groups; the definition is
related to the Poisson formula for the ball $B_{R}:=\left\{  x\in
\mathbb{R}^{n}:\left|  x\right|  <R\right\}  $ which we recall now: Let $R>0$
and $h$ be a function harmonic in the ball $B_{R}$ and continuous on the
closure $\overline{B_{R}};$ then for any $x\in\mathbb{R}^{n}$ with $\left|
x\right|  <R$
\begin{equation}
h\left(  x\right)  =\frac{1}{\omega_{n}}\int_{\mathbb{S}^{n-1}}\frac{\left(
R^{2}-\left|  x\right|  ^{2}\right)  R^{n-2}}{r\left(  R\theta-x\right)  ^{n}%
}h\left(  R\theta\right)  d\theta, \label{ppoisson}%
\end{equation}
where $\omega_{n}$ denotes the area of $\mathbb{S}^{n-1},$ $\theta
\in\mathbb{S}^{n-1},$ $y=R\theta$, and $r\left(  x\right)  $ is the euclidean
norm of $x$. Note that for fixed $x$ with $\left|  x\right|  <R$ the function
$\rho\longmapsto r\left(  \rho\theta-x\right)  $ defined for $\rho
\in\mathbb{R}$ with $\left|  \rho\right|  >R$ has an analytic continuation for
$\zeta\in\mathbb{C}$ with $\left|  \zeta\right|  >R,$ so we can write
$r\left(  \zeta\theta-x\right)  $ for $\zeta\in\mathbb{C}$ with $\left|
\zeta\right|  >R.$ The following \emph{Cauchy type integral formula}, proved
in \cite[p. 125]{ACL83}, is important for our approach: for any polynomial
$u\left(  x\right)  $ and for any $\left|  x\right|  <R$ the following
identity holds
\begin{equation}
u\left(  x\right)  =\frac{1}{2\pi i\omega_{n}}\int_{\Gamma_{R}}\int
_{\mathbb{S}^{n-1}}\frac{\zeta^{n-1}}{r\left(  \zeta\theta-x\right)  ^{n}%
}u\left(  \zeta\theta\right)  d\theta d\zeta\label{holom}%
\end{equation}
where the contour $\Gamma_{R}\left(  t\right)  =R\cdot e^{it}$ for
$t\in\left[  0,2\pi\right]  $. A similar result is also valid for holomorphic
functions $u$ defined on the so-called harmonicity hull of $B_{R}$; since we
need (\ref{holom}) only for polynomials we refer the reader to \cite[p.
125]{ACL83} for details.

Assume now that $\mu$ is a measure with support in the closed ball $\left\{
x\in\mathbb{R}^{n}:\left|  x\right|  \leq R\right\}  .$ The \emph{multivariate
Markov transform} $\widehat{\mu}$ of $\mu$ is a function defined for all
$\theta\in\mathbb{S}^{n-1}$ and all $\zeta\in\mathbb{C}$ with $\left|
\zeta\right|  >R$ by the formula
\begin{equation}
\widehat{\mu}\left(  \zeta,\theta\right)  =\frac{1}{\omega_{n}}\int
_{\mathbb{R}^{n}}\frac{\zeta^{n-1}}{r\left(  \zeta\theta-x\right)  ^{n}}%
d\mu\left(  x\right)  . \label{defmst}%
\end{equation}
Since $\zeta\mapsto r\left(  \zeta\theta-x\right)  $ has no zeros for $\left|
\zeta\right|  >R$ the function $\zeta\mapsto\widehat{\mu}\left(  \zeta
,\theta\right)  $ is defined for all $\left|  \zeta\right|  >R.$ In the first
Section we shall show that the multivariate Markov transform $\widehat{\mu}$
determines the measure $\mu$ uniquely, cf. Theorem \ref{T2}.

Our second main innovation is the introduction of the notion of the
\emph{function }$Q_{P}\left(  \zeta,\theta\right)  $ \emph{of the second kind}
with respect to a given polynomial $P\left(  x\right)  $ which is the
multivariate analogue of (\ref{secondkind}), defined by
\begin{equation}
Q_{P}\left(  \zeta,\theta\right)  =\int_{\mathbb{R}^{n}}\frac{P\left(
\zeta\theta\right)  -P\left(  x\right)  }{r\left(  \zeta\theta-x\right)  ^{n}%
}\zeta^{n-1}d\mu\left(  x\right)  \text{ } \label{secondkind2}%
\end{equation}
for all $\left|  \zeta\right|  >R,\theta\in\mathbb{S}^{n-1}.$ Let us emphasize
that $Q_{P}$ is in general \emph{not} a polynomial. However, we shall show the
surprising and interesting result that the function $r\theta\mapsto
r^{-\left(  n-1\right)  }Q_{P}\left(  r\theta\right)  $ is a
\emph{polyharmonic} function of order $\leq\deg P\left(  x\right)  $ where
$\deg$ denotes the usual total degree of a polynomial.

One further $\emph{main}$ \emph{result }of the paper, Theorem \ref{T20}, is
concerned with measures $\mu$ having their supports in algebraic sets: Let us
assume that the measure $\mu$ has support in $K_{P}\left(  R\right)  .$ Then
the Markov transform $\widehat{\mu}$ has the representation
\begin{equation}
\widehat{\mu}\left(  \zeta,\theta\right)  =\frac{Q_{P}\left(  \zeta
,\theta\right)  }{P\left(  \zeta\theta\right)  }\qquad\text{for }\left|
\zeta\right|  >R, \label{QdividedP}%
\end{equation}
where $Q_{P}$ is the function of second kind with respect to $P\left(
x\right)  .$ The reverse statement holds as well, i.e. if the measure $\mu$
with $\operatorname*{supp}\left(  \mu\right)  \subset\overline{B_{R}}$
satisfies (\ref{QdividedP}) for some polynomial $P$ where $Q_{P}$ is defined
by (\ref{secondkind2})$,$ then $\operatorname*{supp}\left(  \mu\right)
\subset K_{P}\left(  R\right)  .$ By means of these characterizations we can
deduce our main result Theorem \ref{ThmMain}.

\section{The multivariate Markov transform}

Recall that the univariate Markov transform has, for $\left|  \zeta\right|
>R,$ the asymptotic expansion
\begin{equation}
\widehat{\sigma}\left(  \zeta\right)  =\sum_{k=0}^{\infty}\frac{1}{\zeta
^{k+1}}\int_{-\infty}^{\infty}t^{k}d\sigma\left(  t\right)  .
\label{sasymptotic}%
\end{equation}
Let $\Gamma_{R}$ denote the contour in $\mathbb{C}$ defined by $\Gamma
_{R}\left(  t\right)  =R\cdot e^{it}$ for $t\in\left[  0,2\pi\right]  $. By
means of standard facts from complex analysis the following identity may be
proved,
\begin{equation}
M\left(  p\right)  :=\frac{1}{2\pi i}\int_{\Gamma_{R_{1}}}p\left(
\zeta\right)  \widehat{\sigma}\left(  \zeta\right)  d\zeta=\int_{-R}%
^{R}p\left(  x\right)  d\sigma\left(  x\right)  \label{Mintegral}%
\end{equation}
for all polynomials $p$ and any $R_{1}>R$.

In this section we want to show that similar results hold for the multivariate
Markov transform $\widehat{\mu};$ in particular the following is the analogue
of formula (\ref{Mintegral}) in the multivariate case:

\begin{proposition}
\label{P1}\label{stielt2}Let $\mu$ be a signed measure over $\mathbb{R}^{n}$
with support in $\overline{B_{R}}$ and let $R_{1}>R.$ Then for every polnomial
$P\left(  x\right)  $
\begin{equation}
M_{\mu}\left(  P\right)  :=\frac{1}{2\pi i}\int_{\Gamma_{R_{1}}}%
\int_{\mathbb{S}^{n-1}}P\left(  \zeta\theta\right)  \widehat{\mu}\left(
\zeta,\theta\right)  d\zeta d\theta=\int_{\mathbb{R}^{n}}P\left(  x\right)
d\mu\left(  x\right)  . \label{MMintegral}%
\end{equation}

\end{proposition}

%

\proof
Replace $\widehat{\mu}\left(  \zeta,\theta\right)  $ in (\ref{MMintegral}) by
(\ref{defmst}) and interchange integration. Then
\begin{equation}
M_{\mu}\left(  P\right)  =\int_{\mathbb{R}^{n}}\frac{1}{2\pi i\omega_{n}}%
\int_{\Gamma_{R_{1}}}\int_{\mathbb{S}^{n-1}}P\left(  \zeta\theta\right)
\frac{\zeta^{n-1}}{r\left(  \zeta\theta-x\right)  ^{n}}d\zeta d\theta
d\mu\left(  x\right)  . \label{MMMintegral}%
\end{equation}
According to (\ref{holom}) we obtain $M_{\mu}\left(  P\right)  =\int P\left(
x\right)  d\mu\left(  x\right)  .$%
\endproof

\begin{theorem}
\label{T1}\label{uniquenessST}\label{T2}Let $\mu,\nu$ be finite signed
measures over $\mathbb{R}^{n}$ with compact support. If the multivariate
Markov transforms of $\mu$ and $\nu$ coincide for large $\zeta$, i.e., if
there exists $R>0$ such that $\widehat{\mu}\left(  \zeta,\theta\right)
=\widehat{\nu}\left(  \zeta,\theta\right)  $ for all $\left|  \zeta\right|  >R
$ and for all $\theta\in\mathbb{S}^{n-1},$ then $\mu$ and $\nu$ are identical.
\end{theorem}

%

\proof
Since the multivariate Markov transforms coincide for large $\left|
\zeta\right|  $ it is clear that the functionals $M_{\mu}$ and $M_{\nu}$ in
(\ref{MMintegral}) are identical by taking the radius $R_{1}$ of the path
$\Gamma_{R_{1}}$ large enough. Then Proposition \ref{P1} shows that $\int
P\left(  x\right)  d\mu\left(  x\right)  =\int P\left(  x\right)  d\nu\left(
x\right)  $ for all polynomials $P\left(  x\right)  .$ Further we apply a
standard argument: since $\mu$ and $\nu$ have compact supports we may apply
the Stone--Weierstrass theorem according to which the polynomials are dense in
the space $C\left(  \operatorname*{supp}\left(  \mu\right)  \cup
\operatorname*{supp}\left(  \nu\right)  \right)  $ which implies by the
Hahn--Banach theorem that $\mu=\nu.$%
\endproof

Next we want to determine the asymptotic expansion of the multivariate Markov
transform and we need some notations from harmonic analysis; for a detailed
account we refer to \cite{ABR92} or \cite{StWe71}. Recall that a function
$Y:\mathbb{S}^{n-1}\rightarrow\mathbb{C}$ is called a \emph{spherical
harmonic} of degree $k\in\mathbb{N}_{0}$ if there exists a \emph{homogeneous}
\emph{harmonic} polynomial $P\left(  x\right)  $ of degree $k$ (in general,
with complex coefficients) such that $P\left(  \theta\right)  =Y\left(
\theta\right)  $ for all $\theta\in\mathbb{S}^{n-1}.$\footnote{One may
restrict the attention to real valued spherical harmonics and this does not
change the results essentially.} Throughout the paper we assume that
$Y_{k,m}\left(  x\right)  $, $m=1,...,a_{k},$ is a basis of the set of all
harmonic homogeneous polynomials of degree $k$ which are orthonormal with
respect to scalar product
\[
\left\langle f,g\right\rangle _{\mathbb{S}^{n-1}}:=\int_{\mathbb{S}^{n-1}%
}f_{m}\left(  \theta\right)  \overline{g\left(  \theta\right)  }d\theta.
\]
For a continuous function $f:\mathbb{S}^{n-1}\rightarrow\mathbb{C}$ we define
the \emph{Laplace-Fourier series} by
\[
f\left(  \theta\right)  =\sum_{k=0}^{\infty}\sum_{m=1}^{a_{k}}f_{k,m}%
Y_{k,m}\left(  \theta\right)
\]
and $f_{k,m}=\int_{\mathbb{S}^{n-1}}f\left(  \theta\right)  \overline
{Y_{k,m}\left(  \theta\right)  }d\theta$ are the \emph{Laplace-Fourier
coefficients} of $f.$

Using the \emph{Gauss decomposition} of a polynomial (see Theorem 5.5 in
\cite{ABR92}) it is easy to see that the system
\[
\left|  x\right|  ^{2t}Y_{k,m}\left(  x\right)  ,t,k\in\mathbb{N}%
_{0},m=1,...,a_{k}%
\]
is a basis of the set of all polynomials. The numbers
\begin{equation}
c_{t,k,m}:=\int_{\mathbb{R}^{n}}\left|  x\right|  ^{2t}\overline
{Y_{k,m}\left(  x\right)  }d\mu\left(  x\right)  ,\quad t,k\in\mathbb{N}%
_{0},m=1,...,a_{k} \label{distributed}%
\end{equation}
are sometimes called the \emph{distributed moments}, see \cite{kounchev87}.
For a treatment and formulation of the \emph{multivariate moment problem} we
refer to \cite{Fugl83}, see also \cite{StSz98}.

\begin{theorem}
\label{T3}Let $\mu$ be a signed measure over $\mathbb{R}^{n}$ with support in
the closed ball $\overline{B_{R}}$. Then for all $\left|  \zeta\right|  >R$
and for all $\theta\in\mathbb{S}^{n-1}$ the following relation holds
\begin{equation}
\widehat{\mu}\left(  \zeta,\theta\right)  =\sum_{t=0}^{\infty}\sum
_{k=0}^{\infty}\sum_{m=1}^{a_{k}}\frac{Y_{k,m}\left(  \theta\right)  }%
{\zeta^{2t+k+1}}\int_{\mathbb{R}^{n}}\left|  x\right|  ^{2t}\overline
{Y_{k,m}\left(  x\right)  }d\mu\left(  x\right)  \label{muhatrep2}%
\end{equation}

\end{theorem}

%

\proof
A zonal harmonic of degree $k$ with pole $\theta\in\mathbb{S}^{n-1}$ is the
unique spherical harmonic $Z_{\theta}^{\left(  k\right)  }$ of degree $k$ such
that for all spherical harmonics $Y$ of degree $k$ the relation $Y\left(
\theta\right)  =\int_{\mathbb{S}^{n-1}}Z_{\theta}^{\left(  k\right)  }\left(
\eta\right)  Y\left(  \eta\right)  d\eta$ holds. Let $p_{n}\left(
\theta,x\right)  =\frac{1}{\omega_{n}}\frac{1-\left|  x\right|  ^{2}}{\left|
x-\theta\right|  ^{n}}$ be the Poisson kernel for $0\leq\left|  x\right|
<1=\left|  \theta\right|  .$ Theorem 2.10 in \cite[p. 145]{StWe71} gives
$p_{n}\left(  \theta,x\right)  =\sum_{k=0}^{\infty}\left|  x\right|
^{k}Z_{\theta}^{\left(  k\right)  }\left(  x^{\prime}\right)  $ for all
$\theta,x^{\prime}\in\mathbb{S}^{n-1},$ where $x=\left|  x\right|  \cdot
x^{\prime},\ \left|  x\right|  <1.$ Lemma 2.8 in \cite{StWe71} shows that
$Z_{\theta}^{\left(  k\right)  }\left(  x^{\prime}\right)  =\sum_{m=1}^{a_{k}%
}\overline{Y_{k,m}\left(  x^{\prime}\right)  }Y_{k,m}\left(  \theta\right)  $
where $x^{\prime},\theta\in\mathbb{S}^{n-1},$ so
\begin{equation}
p_{n}\left(  \theta,x\right)  =\sum_{k=0}^{\infty}\sum_{m=1}^{a_{k}}\left|
x\right|  ^{k}\overline{Y_{k,m}\left(  x^{\prime}\right)  }Y_{k,m}\left(
\theta\right)  . \label{eqpn}%
\end{equation}
for $\left|  x\right|  <1.$ Let $R$ be as in the theorem, and replace now $x$
in (\ref{eqpn}) by $x/\rho$, $\rho\in\mathbb{R}$ such that $\left|  x\right|
<R<\rho;$ one obtains that
\begin{equation}
\frac{1}{\omega_{n}}\frac{\rho^{n-2}\left(  \rho^{2}-\left|  x\right|
^{2}\right)  }{r\left(  \rho\theta-x\right)  ^{n}}=\sum_{k=0}^{\infty}%
\sum_{m=1}^{a_{k}}\frac{1}{\rho^{k}}\overline{Y_{k,m}\left(  x\right)
}Y_{k,m}\left(  \theta\right)  . \label{gleich2}%
\end{equation}
The real variable $\rho$ can now be replaced by a complex variable $\zeta$
with $\left|  \zeta\right|  >R.$ We multiply by $\zeta\left(  \zeta
^{2}-\left|  x\right|  ^{2}\right)  ^{-1},$ and integrate integrate over the
closed ball $\overline{B_{R}}$ with respect to $\mu$. This gives
\begin{equation}
\widehat{\mu}\left(  \zeta,\theta\right)  =\sum_{k=0}^{\infty}\sum
_{m=1}^{a_{k}}Y_{k,m}\left(  \theta\right)  \zeta^{-k+1}\int_{\mathbb{R}^{n}%
}\frac{\overline{Y_{k,m}\left(  x\right)  }}{\zeta^{2}-\left|  x\right|  ^{2}%
}d\mu\left(  x\right)  , \label{gleich4}%
\end{equation}
and we have determined the Laplace-Fourier series of $\theta\longmapsto
\widehat{\mu}\left(  \zeta,\theta\right)  .$ Since $\left|  \zeta\right|
>R\geq\left|  x\right|  $ we can expand $1/(1-\frac{\left|  x\right|  ^{2}%
}{\zeta^{2}})$ in a geometric series and we obtain
\begin{equation}
\widehat{\mu}\left(  \zeta,\theta\right)  =\sum_{k=0}^{\infty}\sum
_{m=1}^{a_{k}}\frac{Y_{k,m}\left(  \theta\right)  }{\zeta^{k+1}}%
\int_{\mathbb{R}^{n}}\overline{Y_{k,m}\left(  x\right)  }\left(  \sum
_{t=0}^{\infty}\frac{\left|  x\right|  ^{2t}}{\zeta^{2t}}\right)  d\mu\left(
x\right)  . \label{gleich5}%
\end{equation}
After interchanging summation and integration the claim is obvious.%
\endproof

\section{The function of the second kind}

In the following we want to a give a multivariate analogue of the polynomial
of second kind. It turns out that in the multivariate case the corresponding
definition does not lead to a polynomial but to a polyharmonic function
$Q_{P}\left(  \zeta,\theta\right)  $ which is defined only for all $\left|
\zeta\right|  >R,\theta\in\mathbb{S}^{n-1}$.

\begin{definition}
Let $P\left(  x\right)  $ be a polynomial and $\mu$ be a non-negative measure
with support in $\overline{B_{R}}.$ Then the \emph{function }$Q_{P}\left(
\zeta,\theta\right)  $ \emph{of the second kind} is defined by
\[
Q_{P}\left(  \zeta,\theta\right)  =\frac{1}{\omega_{n}}\int_{\mathbb{R}^{n}%
}\frac{P\left(  \zeta\theta\right)  -P\left(  x\right)  }{r\left(  \zeta
\theta-x\right)  ^{n}}\zeta^{n-1}d\mu\left(  x\right)
\]
for all $\left|  \zeta\right|  >R,\theta\in\mathbb{S}^{n-1}.$ Similarly we
define the function $R_{P}\left(  \zeta,\theta\right)  $ by
\[
R_{P}\left(  \zeta,\theta\right)  =\frac{1}{\omega_{n}}\int_{\mathbb{R}^{n}%
}\frac{P\left(  x\right)  }{r\left(  \zeta\theta-x\right)  ^{n}}\zeta
^{n-1}d\mu\left(  x\right)
\]
for all $\left|  \zeta\right|  >R,\theta\in\mathbb{S}^{n-1}.$
\end{definition}

The last definitions immediately give the identity $\ $
\begin{equation}
P\left(  \zeta\theta\right)  \widehat{\mu}\left(  \zeta,\theta\right)
=Q_{P}\left(  \zeta,\theta\right)  +R_{P}\left(  \zeta,\theta\right)  .
\label{form1}%
\end{equation}

\begin{theorem}
\label{T6}Let $P\left(  x\right)  $ be a polynomial, $\mu$ be a signed measure
with support in $\overline{B_{R}}$ and $Q_{P}\left(  \zeta,\theta\right)  $
the function of the second kind. Then for any $R_{1}>R$ and for each
polynomial $h\left(  x\right)  $
\begin{equation}
\frac{1}{2\pi i}\int_{\Gamma_{R_{1}}}\int_{\mathbb{S}^{n-1}}h\left(
\zeta\theta\right)  Q_{P}\left(  \zeta,\theta\right)  d\zeta d\theta=0.
\label{IHO}%
\end{equation}

\end{theorem}

%

\proof
Let us denote the integral in (\ref{IHO}) by $I\left(  h\right)  .$ By
(\ref{form1}) we obtain that $I\left(  h\right)  =I_{1}\left(  h\right)
-I_{2}\left(  h\right)  $ where
\begin{align}
I_{1}\left(  h\right)   &  =\frac{1}{2\pi i}\int_{\Gamma_{R_{1}}}%
\int_{\mathbb{S}^{n-1}}h\left(  \zeta\theta\right)  P\left(  \zeta
\theta\right)  \widehat{\mu}\left(  \zeta,\theta\right)  d\zeta d\theta
,\label{zeile1}\\
I_{2}\left(  h\right)   &  =\frac{1}{2\pi i\omega_{n}}\int_{\Gamma_{R_{1}}%
}\int_{\mathbb{S}^{n-1}}h\left(  \zeta\theta\right)  \int_{\mathbb{R}^{n}%
}\frac{P\left(  x\right)  }{r\left(  \zeta\theta-x\right)  ^{n}}\zeta
^{n-1}d\mu\left(  x\right)  d\zeta d\theta. \label{zeile2}%
\end{align}
Proposition \ref{P1} yields $I_{1}\left(  h\right)  =\int_{\mathbb{R}^{n}%
}h\left(  x\right)  P\left(  x\right)  d\mu\left(  x\right)  .$ Change the
integration order in (\ref{zeile2}) and use formula (\ref{holom}). Then we
obtain $I_{2}\left(  h\right)  =I_{1}\left(  h\right)  ,$ therefore $I\left(
h\right)  =0$ which was our claim.%
\endproof

A similar argument as in the proof of formula (\ref{muhatrep2}) proves the following:

\begin{theorem}
The rest function $R_{P}\left(  \zeta,\theta\right)  $ has the asymptotic
expansion
\begin{equation}
\sum_{t=0}^{\infty}\sum_{k=0}^{\infty}\sum_{m=1}^{a_{k}}\frac{Y_{k,m}\left(
\theta\right)  }{\zeta^{2t+k+1}}\int_{\mathbb{R}^{n}}P\left(  x\right)
\left|  x\right|  ^{2t}\overline{Y_{k,m}\left(  x\right)  }d\mu\left(
x\right)  . \label{form3}%
\end{equation}

\end{theorem}

Let us consider now the Laurent series of the function $\zeta\mapsto
R_{P}\left(  \zeta,\theta\right)  $: for $\left|  \zeta\right|  >R,\theta
\in\mathbb{S}^{n-1}$ we can write
\begin{equation}
R_{P}\left(  \zeta,\theta\right)  =\sum_{s=0}^{\infty}r_{s}\left[  P\right]
\left(  \theta\right)  \frac{1}{\zeta^{s+1}}. \label{Rest}%
\end{equation}
From (\ref{form3}), by putting $s=2t+k,$ it follows that
\begin{equation}
r_{s}\left[  P\right]  \left(  \theta\right)  =\sum_{t=0}^{\left[  s/2\right]
}\sum_{m=1}^{a_{s-2t}}Y_{s-2t,m}\left(  \theta\right)  \int_{\mathbb{R}^{n}%
}P\left(  x\right)  \left|  x\right|  ^{2t}\overline{Y_{s-2t,m}\left(
x\right)  }d\mu\left(  x\right)  . \label{neu1}%
\end{equation}
Hence the coefficient function $r_{s}\left(  P\right)  $ is a sum of spherical
harmonics with degree $\leq s$.

We can now formulate a characterization of orthogonality in asymptotic analysis:

\begin{theorem}
\label{Torthog}\label{TO1}Let $\mu$ be a signed measure with compact support
and $P\left(  x\right)  $ be a polynomial. Then $P$ is orthogonal to all
polynomials of degree $<M$ with respect to $\mu$ if and only if
\[
r_{0}\left[  P\right]  =...=r_{M-1}\left[  P\right]  =0
\]
where $r_{s}\left[  P\right]  $ are the functions defined in (\ref{Rest}%
)--(\ref{neu1}).
\end{theorem}

%

\proof
From (\ref{neu1}) we see that $r_{0}(P)=...=r_{M-1}(P)=0$ if and only for all
$s=0,...,M-1$
\[
\int_{\mathbb{R}^{n}}P\left(  x\right)  \left|  x\right|  ^{2t}\overline
{Y_{s-2t,m}\left(  x\right)  }d\mu\left(  x\right)  =0.
\]
But the polynomials $\left|  x\right|  ^{2t}Y_{s-2t,m}\left(  x\right)  $ with
$s=0,...,M-1,$ $t=0,...,\left[  s/2\right]  ,m=1,...,a_{s-2t}$, span up the
space of polynomials of degree $\leq M-1.$%
\endproof

The next theorem, interesting in its own right, is not needed later, and
therefore the proof will be omitted.

\begin{theorem}
\label{T19}Let $\mu$ be a signed measure with compact support and let
$P\left(  x\right)  $ be a polynomial of degree $2N.$ If $P$ is orthogonal to
all polynomials of degree $\leq2N$ and polyharmonic degree $<N$ then
$r_{0}(P)=...=r_{2N-1}(P)=0$ and $r_{2N}\left(  \theta\right)  $ is constant.
\end{theorem}

\section{Polyharmonicity of the function of second kind}

In this Section we want to show that the function $Q_{P}\left(  \zeta
,\theta\right)  $ of the second kind, multiplied by $\zeta^{-\left(
n-1\right)  },$ is a polyharmonic function.

Recall that we have defined $N_{P}=\sup\left\{  d\left(  P\cdot h\right)
:h\text{ harmonic polynomial}\right\}  $ for a polynomial $P\left(  x\right)
.$ In the Appendix we will show that $N_{P}\leq\deg P\left(  x\right)  $ and
an explicit determination of $N_{P}$ will be given there as well.

\begin{proposition}
Let $Y_{k,m},m=1,...,a_{k},$ be an orthonormal basis of the space of all
homogeneous harmonic polynomials. Then
\begin{equation}
N_{P}:=\sup_{k\in\mathbb{N}_{0},m=1,...,a_{k}}d\left(  P\left(  x\right)
Y_{k,m}\left(  x\right)  \right)  . \label{NP}%
\end{equation}

\end{proposition}

%

\proof
Let us denote the right hand side by $M_{P}.$ Then the inequality $M_{P}\leq
N_{P}$ is trivial. For the converse let $h\left(  x\right)  $ be a harmonic
polynomial and write $h\left(  x\right)  =\sum_{k=0}^{N}\sum_{m=1}^{a_{k}%
}\lambda_{k,m}Y_{k,m}\left(  x\right)  .$ Then
\[
d\left(  P\cdot h\right)  \leq\sup_{k\in\mathbb{N}_{0},m=1,...,a_{k}}d\left(
P\left(  x\right)  Y_{k,m}\left(  x\right)  \right)  \ \leq M_{P}.
\]%
\endproof

Note that $N_{P}=\sup_{k\in\mathbb{N}_{0},m=1,...,a_{k}}d\left(  P\left(
x\right)  \overline{Y_{k,m}\left(  x\right)  }\right)  $ since $\overline
{Y_{k,m}},$ $m=1,...,a_{k}$ is an orthonormal basis as well. Now we determine
the asymptotic expansion of the function of the second kind:

\begin{theorem}
\label{Tdecompositionmu}Let $P\left(  x\right)  $ be a polynomial and $\mu$ be
a signed measure with support in $\overline{B_{R}.}$ Then $\theta\mapsto
Q_{P}\left(  \zeta,\theta\right)  ,$ the function of the second kind,
possesses a Laplace-Fourier series of the form
\begin{equation}
Q_{P}\left(  \zeta,\theta\right)  =\sum_{k=0}^{\infty}\sum_{m=1}^{a_{k}}%
\frac{1}{\zeta^{k-1}}p_{k,m}\left(  \zeta^{2}\right)  Y_{k,m}\left(
\theta\right)  \label{EP}%
\end{equation}
where $p_{k,m}\left(  t\right)  $ are univariate polynomials of degree
strictly smaller than $N_{k,m}:=d\left(  P\left(  x\right)  Y_{k,m}\left(
x\right)  \right)  .$ The function $Q_{P}\left(  \zeta,\theta\right)  $ of the
second kind depends on those distributed moments
\begin{equation}
\int_{\mathbb{R}^{n}}h\left(  x\right)  \left|  x\right|  ^{2t}d\mu\left(
x\right)  \label{dist1}%
\end{equation}
where $t\leq\sup_{k\in\mathbb{N}_{0}}\deg p_{k,m}$ and $h\left(  x\right)  $
is a harmonic polynomial.
\end{theorem}

%

\proof
For each fixed $\zeta$ with $\left|  \zeta\right|  >R$ the function
$\theta\mapsto Q_{P}\left(  \zeta,\theta\right)  $ possesses a Laplace-Fourier
expansion, say
\[
Q_{P}\left(  \zeta,\theta\right)  =\sum_{k=0}^{\infty}\sum_{m=1}^{a_{k}}%
e_{km}\left(  \zeta\right)  Y_{k,m}\left(  \theta\right)
\]
Recall that $Q_{P}\left(  \zeta,\theta\right)  =$ $P\left(  \zeta
\theta\right)  \widehat{\mu}\left(  \zeta,\theta\right)  -R_{P}\left(
\zeta,\theta\right)  .$ Formula (\ref{form3}) yields the Laplace-Fourier
expansion of $\theta\mapsto R_{P}\left(  \zeta,\theta\right)  :$ in
(\ref{form3}) one computes the sum over the variable $t$ obtaining
\begin{equation}
R_{P}\left(  \zeta,\theta\right)  =\sum_{k=0}^{\infty}\sum_{m=1}^{a_{k}%
}Y_{k,m}\left(  \theta\right)  \frac{1}{\zeta^{k-1}}\int_{\mathbb{R}^{n}}%
\frac{P\left(  x\right)  \overline{Y_{k,m}\left(  x\right)  }}{\zeta
^{2}-\left|  x\right|  ^{2}}d\mu\left(  x\right)  . \label{Laplace1}%
\end{equation}
The Laplace-Fourier coefficients of $\theta\mapsto P\left(  \zeta
\theta\right)  \widehat{\mu}\left(  \zeta,\theta\right)  $ are given through
\begin{equation}
f_{k,m}\left(  \zeta\right)  :=\int_{\mathbb{S}^{n-1}}P\left(  \zeta
\theta\right)  \widehat{\mu}\left(  \zeta,\theta\right)  \overline
{Y_{k,m}\left(  \theta\right)  }d\theta. \label{fkm}%
\end{equation}
Let us write $P\left(  x\right)  \overline{Y_{k,m}\left(  x\right)  }$ in the
Gau\ss \ decomposition, see Theorem $5.5$ in \cite{ABR92}, in the form
\begin{equation}
P\left(  x\right)  \overline{Y_{k,m}\left(  x\right)  }=\sum_{j=0}^{N_{k,m}%
}h_{j,k,m}\left(  x\right)  \left|  x\right|  ^{2j}, \label{PYgauss}%
\end{equation}
where $h_{j,k,m}$ are harmonic polynomials and $N_{k,m}$ is the polyharmonic
degree of $P\left(  x\right)  Y_{k,m}\left(  x\right)  .$ Then (\ref{fkm}) and
(\ref{PYgauss}) yield
\begin{align*}
f_{k,m}\left(  \zeta\right)   &  =\frac{1}{\zeta^{k}}\int_{\mathbb{S}^{n-1}%
}P\left(  \zeta\theta\right)  \zeta^{k}\overline{Y_{k,m}\left(  \theta\right)
}\widehat{\mu}\left(  \zeta,\theta\right)  d\theta\\
&  =\frac{1}{\zeta^{k}}\sum_{j=0}^{N_{k,m}}\zeta^{2j}\int_{\mathbb{S}^{n-1}%
}h_{j,k,m}\left(  \zeta\theta\right)  \widehat{\mu}\left(  \zeta
,\theta\right)  d\theta\\
&  =\frac{1}{\zeta^{k}}\sum_{j=0}^{N_{k,m}}\zeta^{2j}\int_{\mathbb{R}^{n}}%
\int_{\mathbb{S}^{n-1}}h_{j,k,m}\left(  \zeta\theta\right)  \frac{1}%
{\omega_{n}}\frac{\zeta^{n-1}}{r\left(  \zeta\theta-x\right)  ^{n}}d\theta
d\mu\left(  x\right)  .
\end{align*}
Since $h_{j,k,m}$ is a harmonic polynomial the Poisson formula shows that for
real $\zeta>R$ holds
\[
h_{j,k,m}\left(  x\right)  =\frac{1}{\omega_{n}}\int_{\mathbb{S}^{n-1}%
}h_{j,k,m}\left(  \zeta\theta\right)  \frac{\zeta^{n-2}\left(  \zeta
^{2}-\left|  x\right|  ^{2}\right)  }{r\left(  \zeta\theta-x\right)  ^{n}%
}d\theta.
\]
Since the integrand is holomorphic in $\zeta$ this holds for all complex
values $\zeta$ with $\left|  \zeta\right|  >R$ as well. Thus
\begin{equation}
f_{k,m}\left(  \zeta\right)  =\frac{1}{\zeta^{k}}\sum_{j=0}^{N_{k,m}}%
\zeta^{2j}\int_{\mathbb{R}^{n}}\frac{\zeta}{\zeta^{2}-\left|  x\right|  ^{2}%
}h_{j,k,m}\left(  x\right)  d\mu\left(  x\right)  \label{egunbekannt}%
\end{equation}
are the Laplace Fourier coefficients of $\theta\mapsto P\left(  \zeta
\theta\right)  \widehat{\mu}\left(  \zeta,\theta\right)  .$

Replace now $P\left(  x\right)  \overline{Y_{k,m}\left(  x\right)  }$ in
(\ref{Laplace1}) by the right hand side of (\ref{PYgauss}) and take the
difference of the Laplace-Fourier coefficients we computed so far. Then the
Laplace-Fourier coefficients of $Q_{P}\left(  \zeta,\theta\right)  $ are given
by
\[
e_{k,m}\left(  \zeta\right)  =\frac{1}{\zeta^{k-1}}\sum_{j=0}^{N_{k,m}}%
\int_{\mathbb{R}^{n}}\frac{1}{\zeta^{2}-\left|  x\right|  ^{2}}h_{j,k,m}%
\left(  x\right)  \left(  \zeta^{2j}-\left|  x\right|  ^{2j}\right)
d\mu\left(  x\right)  .
\]
Note that for $j=0$ the summand ist just zero. For $j\geq1$ we have
\[
\frac{\zeta^{2j}-\left|  x\right|  ^{2j}}{\zeta^{2}-\left|  x\right|  ^{2}%
}=\left|  x\right|  ^{2\left(  j-1\right)  }+\left|  x\right|  ^{2\left(
j-1\right)  }\zeta^{2}+...+\zeta^{2\left(  j-1\right)  }.
\]
We conclude that $\zeta\mapsto\zeta^{k-1}e_{k,m}\left(  \zeta\right)
=:P_{k,m}\left(  \zeta^{2}\right)  $ is a polynomial in $\zeta^{2}$ of degree
at most $N_{k,m}-1.$ It follows that $e_{k,m}\left(  \zeta\right)  $ can be
computed if we know all moments of the form (\ref{dist1}) where $t\leq\deg
p_{k,m}$ and $h\left(  x\right)  $ is a harmonic polynomial. The proof is
complete.%
\endproof

From this we have the following interesting consequence

\begin{corollary}
Let $P\left(  x\right)  $ be a polynomial, $\mu$ be a signed measure with
support in $\overline{B_{R}}$ and $Q_{P}\left(  \zeta,\theta\right)  $ be the
corresponding function of the second kind. Then the function $r\theta\mapsto
r^{-\left(  n-1\right)  }Q_{P}\left(  r\theta\right)  $ defined for $r>R$ and
$\theta\in\mathbb{S}^{n-1},$ is a polyharmonic function of polyharmonic degree
$<N_{P}$ where $N_{P}$ is defined in (\ref{NP}).
\end{corollary}

%

\proof
By the last theorem the function $\theta\mapsto r^{-\left(  n-1\right)  }%
Q_{P}\left(  r\theta\right)  $ has the following Laplace-Fourier expansion
\[
f\left(  r\theta\right)  :=r^{-\left(  n-1\right)  }Q_{P}\left(
r\theta\right)  =\sum_{k=0}^{\infty}\sum_{m=1}^{a_{k}}\frac{1}{r^{n+k-2}%
}p_{k,m}\left(  r^{2}\right)  Y_{k,m}\left(  \theta\right)
\]
Let us define the differential operator
\begin{equation}
L_{\left(  k\right)  }:=\frac{d^{2}}{dr^{2}}+\frac{n-1}{r}\frac{d}{dr}%
-\frac{k\left(  k+n-2\right)  }{r^{2}}. \label{LLgleichLL}%
\end{equation}
It is known that a function $g\left(  r\theta\right)  $ is a solution of
$\Delta^{p}g\left(  x\right)  =0$ if and only if the coefficient functions
$g_{k,m}\left(  r\right)  $ of its Laplace-Fourier expansion are solutions of
the equation $\left[  L_{\left(  k\right)  }\right]  ^{p}g_{k,m}\left(
r\right)  =0;$ an elaboration of these classical results can be found in
\cite{Koun00}. Further the polynomials $r^{j}$ with
$j=-k-n+2,-k-n+4,...,-k-n+2p$ are solutions of this equation. It follows that
\[
f_{k,m}\left(  r\right)  =\frac{1}{r^{n+k-2}}p_{k,m}\left(  r^{2}\right)
\]
are solutions of the equation $\left[  L_{\left(  k\right)  }\right]
^{p}g_{k,m}\left(  r\right)  =0$ when $p\geq N_{k}.$ The proof is complete.%
\endproof

\section{Measures with algebraic support}

A measure $\mu$ over $\mathbb{R}^{n}$ is \emph{algebraically supported} if the
support of the measure is contained in an algebraic set, i.e. if the support
of $\mu$ is contained in $P^{-1}\left(  0\right)  $ for some polynomial
$P\left(  x\right)  .$ This is equivalent to the statement that $\int P^{\ast
}P\left(  x\right)  d\mu\left(  x\right)  =0$ where $P^{\ast}\left(  x\right)
:=\overline{P\left(  x\right)  }$ for $x\in\mathbb{R}^{n}.$ The Cauchy-Schwarz
inequality implies that
\[
\left|  \int PQd\mu\right|  ^{2}\leq\int PP^{\ast}d\mu\cdot\int Q^{\ast}%
Qd\mu=0.
\]
It follows that $P$ is orthogonal to \emph{all} polynomials $Q$ with respect
to $\mu.$

In the one-dimensional case a measure $\mu$ has algebraic support if and only
if the support is finite. Further this is equivalent to the property that the
Markov transform is a rational function. As we shall see, in the multivariate
case all these properties will be different.

\begin{theorem}
\label{T20}Let $\mu$ be a measure with support in $\overline{B_{R}}$ and let
$P\left(  x\right)  $ be a polynomial. Then$\ \mu$ has support in
$P^{-1}\left\{  0\right\}  $ if and only if
\begin{equation}
P\left(  \zeta\theta\right)  \widehat{\mu}\left(  \zeta,\theta\right)
=Q_{P}\left(  \zeta,\theta\right)  \text{ for all }\theta\in\mathbb{S}%
^{n-1},\left|  \zeta\right|  >R, \label{egleich}%
\end{equation}
where $Q_{P}\left(  \zeta,\theta\right)  $ is the function of the second kind.
\end{theorem}

%

\proof
If $\mu$ has support in $P^{-1}\left\{  0\right\}  $ it follows that the rest
function $R_{P}\left(  \zeta,\theta\right)  $ is equal to zero and
(\ref{egleich}) is evident. For the converse assume that $P\left(  \zeta
\theta\right)  \widehat{\mu}\left(  \zeta,\theta\right)  =Q_{P}\left(
\zeta,\theta\right)  .$ By Proposition \ref{P1} and Theorem \ref{T6}
\begin{align*}
\int P^{\ast}Pd\mu &  =\frac{1}{2\pi i}\int_{\Gamma_{R_{1}}}\int
_{\mathbb{S}^{n-1}}P^{\ast}\left(  \zeta\theta\right)  P\left(  \zeta
\theta\right)  \widehat{\mu}\left(  \zeta,\theta\right)  d\zeta d\theta\\
&  =\frac{1}{2\pi i}\int_{\Gamma_{R_{1}}}\int_{\mathbb{S}^{n-1}}P^{\ast
}\left(  \zeta\theta\right)  Q_{P}\left(  \zeta,\theta\right)  d\zeta
d\theta=0.
\end{align*}
It follows that $\mu$ has support in $P^{-1}\left\{  0\right\}  .$%
\endproof

The same proof shows that $\int P^{\ast}Pd\mu=0$ if we know that for each
fixed $\theta$ the map $\zeta\mapsto P\left(  \zeta\theta\right)  \widehat
{\mu}\left(  \zeta,\theta\right)  $ is a polynomial in the variable $\zeta$
(since the integral over $\Gamma_{R_{1}}$ is already zero). Hence we have
proved that for a measure $\mu$ with compact support the following implication
holds
\[
\quad\zeta\widehat{\mu}\left(  \zeta,\theta\right)  \text{ rational}%
\quad\Rightarrow\quad\text{supp}\left(  \mu\right)  \text{ is contained in an
algebraic set,}%
\]
where rationality of $\widehat{\mu}\left(  \zeta,\theta\right)  $ means that
it is a quotient of two polynomial $Q\left(  x\right)  $ and $P\left(
x\right)  .$ Not very surprisingly, the converse is not true as the following
result shows (where we choose for example $\sigma$ to be equal to the Lebesgue
measure on the unit interval):

\begin{proposition}
Let $\sigma$ be a measure $\sigma$ over $\mathbb{R}$ with compact support,
$\delta_{0}$ the Dirac measure over $\mathbb{R}$ at the point $0$ and let
$\mu=\sigma\otimes\delta_{0}.$ Then the multivariate Markov transform is given
by
\begin{equation}
\widehat{\sigma\otimes\delta_{0}}\left(  \zeta,e^{it}\right)  =\frac{1}%
{\omega_{2}}\sum_{l=0}^{\infty}\int x^{l}d\sigma\left(  x\right)  \frac
{\sin\left(  l+1\right)  t}{\sin t}\frac{1}{\zeta^{l+1}}. \label{sigdelt}%
\end{equation}
Then $\mu$ has algebraic support but its multivariate Markov transform
$\widehat{\sigma\otimes\delta_{0}}$ is rational if and only if the measure
$\sigma$ has finite support.
\end{proposition}

%

\proof
Let $\theta=e^{it}$ with $t\in\mathbb{R}.$ It is straightforward to verfy
that
\begin{align*}
\widehat{\sigma\otimes\delta_{0}}\left(  \zeta,\theta\right)   &  =\frac
{1}{\omega_{2}}\int_{\mathbb{R}^{2}}\frac{\zeta}{r\left(  \zeta\theta-\left(
x,y\right)  \right)  ^{2}}d\left(  \sigma\otimes\delta_{0}\right) \\
&  =\frac{1}{\omega_{2}}\int_{-\infty}^{\infty}\frac{\zeta}{\zeta^{2}-2\zeta
x\cos t+x^{2}}d\sigma.
\end{align*}
Note that
\[
\frac{2i\zeta\sin t}{\zeta^{2}-2\zeta x\cos t+x^{2}}=\frac{1}{\zeta
\overline{\theta}-x}-\frac{1}{\zeta\theta-x}.
\]
Define for the measure $\sigma$ the one-dimensional Markov transform by
$\widetilde{\sigma}\left(  \zeta\right)  =\int\frac{1}{\zeta-x}d\sigma\left(
x\right)  .$ Then $2i\omega_{2}\sin t\cdot\widehat{\sigma\otimes\delta_{0}%
}\left(  \zeta,\theta\right)  =\widetilde{\sigma}\left(  \zeta\overline
{\theta}\right)  -\widetilde{\sigma}\left(  \zeta\theta\right)  $ and the
asymptotic expansion of $\widetilde{\sigma}$ leads to (\ref{sigdelt}).

Assume now that $\widehat{\sigma\otimes\delta_{0}}\left(  \zeta,\theta\right)
$ is rational. Then for $t=\pi/2$ the function $\zeta\mapsto\widehat
{\sigma\otimes\delta_{0}}\left(  \zeta,\theta\right)  $ is rational, i.e. that
$\ f\left(  \zeta\right)  :=\sum_{k=0}^{\infty}\int x^{2k}d\mu\left(
x\right)  \frac{1}{\zeta^{2k+1}}$ is a rational function. From the univariate
results it follows that $\mu$ must have finite support.
\endproof

If $\mu$ is a measure with finite support and the dimension $n$ is even then
it is easy to see that $\zeta\widehat{\mu}\left(  \zeta,\theta\right)  $ is a
rational function. The following example shows that the converse is not true:

\begin{example}
Let $\mu$ be the Lebesgue measure on the unit circle $\mathbb{S}^{1}.$ Since
the measure is rotation-invariant it follows that $\widehat{\mu}\left(
\zeta,\theta\right)  =\frac{\zeta}{\zeta^{2}-1}.$ Hence the multivariate
Markov transform $\zeta\widehat{\mu}\left(  \zeta,\theta\right)  $ is a
rational function but $\mu$ is not discrete.
\end{example}

\section{Proof of Theorem 1}%

\proof
In Theorem \ref{Tdecompositionmu} we have seen that $Q_{\mu,P}$ and $Q_{\nu
,P}$ only depends on the moments $c_{t,k,m}$ where $t<N_{P}.$ It follows that
$Q_{\mu,P}=Q_{\nu,P}.$ By Theorem \ref{T20} $P\left(  \zeta\theta\right)
\widehat{\mu}\left(  \zeta,\theta\right)  =Q_{\mu,P}\left(  \zeta
,\theta\right)  $ and $P\left(  \zeta\theta\right)  \widehat{\nu}\left(
\zeta,\theta\right)  =Q_{\nu,P}\left(  \zeta,\theta\right)  $ for all large
$\zeta$ and for all $\theta\in\mathbb{S}^{n-1},$ therefore $P\left(
\zeta\theta\right)  \widehat{\mu}\left(  \zeta,\theta\right)  =P\left(
\zeta\theta\right)  \widehat{\nu}\left(  \zeta,\theta\right)  .$ We want to
conlcude that $\widehat{\mu}\left(  \zeta,\theta\right)  =\widehat{\nu}\left(
\zeta,\theta\right)  $; in that case Theorem \ref{T2} yields $\mu=\nu.$ If
$P\left(  \zeta\theta\right)  $ has no zeros for large $\zeta$ it is clear
that $\widehat{\mu}\left(  \zeta,\theta\right)  =\widehat{\nu}\left(
\zeta,\theta\right)  .$ In the general case, it suffices to show that
$A:=\left\{  \left(  \zeta,\theta\right)  \in\mathbb{C}\times\mathbb{S}%
^{n-1}:P\left(  \zeta\theta\right)  =0\right\}  $ is nowhere dense since then
a continuity argument leads to $\widehat{\mu}\left(  \zeta,\theta\right)
=\widehat{\nu}\left(  \zeta,\theta\right)  .$ This fact will be proven in the
next Proposition.%
\endproof

Just for completeness sake we include the following

\begin{proposition}
The set $A:=\left\{  \left(  \zeta,\theta\right)  \in\mathbb{C}\times
\mathbb{S}^{n-1}:P\left(  \zeta\theta\right)  =0\right\}  $ is closed and has
no interior point, i.e. $A$ is nowhere dense in $\mathbb{C}\times
\mathbb{S}^{n-1}.$
\end{proposition}

%

\proof
Clearly $A$ is closed. Suppose that there $\theta_{0}\in\mathbb{S}^{n-1}$ and
$\zeta_{0}$ such that $P\left(  \zeta\theta\right)  =0$ for all $\zeta$ in a
neighborhood $U$ of $\zeta_{0}$ and for all $\theta$ in a neighborhood $V$ of
$\theta_{0}.$ For fixed $\theta\in V$ it follows that $\zeta\rightarrow
P\left(  \zeta\theta\right)  $ must be the zero polynomial since for all
$\zeta\in U$ (hence uncountably many $\zeta)$ we have $P\left(  \zeta
\theta\right)  =0.$ It follows that $P\left(  \zeta\theta\right)  =0$ for all
$\zeta\in\mathbb{C}$ and for all $\theta\in V.$ Hence $P\left(  x\right)  =0$
for all $x$ in an open set $W$ of $\mathbb{R}^{n}$ and we conclude that $P=0.$%
\endproof

\begin{corollary}
\label{Cdensity}Let $P\left(  x\right)  $ be a polynomial and $N_{P}$ be given
by (\ref{NP}). Then the space
\[
U_{N_{P}}:=\left\{  Q\in\mathcal{P}_{n}:\ \Delta^{N_{P}}Q=0\right\}
\]
is dense in the space $C\left(  K_{P}\left(  R\right)  ,\mathbb{C}\right)  $
of all continuous complex-valued functions on $K_{P}\left(  R\right)  $
endowed with the supremum norm.
\end{corollary}

%

\proof
Since $U_{N_{P}}$ is closed under complex conjugation we may reduce the
problem to the case of real-valued continuous functions. Suppose that
$U_{N_{P}}$ is not dense in $C\left(  K_{P}\left(  R\right)  ,\mathbb{R}%
\right)  . $ By the Hahn-Banach theorem there exists a continuous non-trivial
real-valued functional $L$ which vanishes on $U_{N_{P}}$. By Riesz's Theorem
there exists a signed measures $\sigma$ representing the functional $L$ with
support in $K_{P}.$ By Theorem \ref{ThmMain} we conclude that $\sigma=0,$ a
contradiction.
\endproof

\section{Appendix: The Polyharmonic degree}

We want to list some of the properties of the polyharmonic degree map. Note
that the inequality $d\left(  P+Q\right)  \leq\max\left\{  d\left(  P\right)
,d\left(  Q\right)  \right\}  $ is trivial. In \cite{ACL83} the important
equality
\begin{equation}
d\left(  Q\cdot\left|  x\right|  ^{2}\right)  =d\left(  Q\right)  +d\left(
\left|  x\right|  ^{2}\right)  =d\left(  Q\right)  +1. \label{additiv}%
\end{equation}
is proved for any polyharmonic function defined on a domain containing zero.
The following inequality is implicitly contained in \cite[Theorem 1.2, p.
31]{ACL83}. For completeness we give the short proof.

\begin{proposition}
\label{P23}\label{ppolydeg}Let $f,g$ be harmonic polynomials. Then $d\left(
ff^{\ast}\right)  =\deg f$ and $d\left(  fg\right)  \leq\min\left\{  \deg
f,\deg g\right\}  $
\end{proposition}

%

\proof
Let $\nabla f$ be the gradient of $f.$ Then $\Delta(fg)=(\Delta f)g+2<\nabla
f,\nabla g>+f\Delta g$. If $h$ and $g$ are harmonic it is easy to show by
induction that
\[
\Delta^{p}\left(  fg\right)  =2^{p}\sum_{i_{1},...,i_{p}=1}^{n}(\frac
{\partial}{\partial x_{i_{1}}}...\frac{\partial}{\partial x_{i_{p}}}%
f)(\frac{\partial}{\partial x_{i_{1}}}...\frac{\partial}{\partial x_{i_{p}}%
}g).
\]
Suppose that $s:=\deg f\leq\deg g.$ Then $\frac{\partial^{\beta}}{\partial
x^{\beta}}f=0$ for all $\beta\in\mathbb{N}_{0}^{n}$ with $\left|
\beta\right|  =s+1.$ It follows from the above formula that $\Delta
^{s+1}(fg)=0.$ Hence $d\left(  fg\right)  =s.$ For the first statement note
that by the above $d(ff^{\ast})\leq\deg f.$ Suppose that $\Delta
^{p+1}(ff^{\ast})=0$ for some $p\in\mathbb{N}.$ Then $\sum_{i_{1}%
,...,i_{p+1}=1}^{n}\left|  \frac{\partial}{\partial x_{i_{1}}}...\frac
{\partial}{\partial x_{i_{p+1}}}f\right|  ^{2}=0.$ It follows that
$\frac{\partial^{\beta}}{\partial x^{\beta}}f=0$ for all $\beta\in
\mathbb{N}_{0}^{n}$ with $\left|  \beta\right|  =p+1.$ Hence $\deg f\leq p$
and we have proved that $\deg f\leq d\left(  ff^{\ast}\right)  .$ The proof is
complete.%
\endproof

Now we can prove the following

\begin{corollary}
Let $Y_{k}$ be a harmonic homogeneous polynomial of degree $k$ and $P\left(
x\right)  $ be a polynomial with the Gau\ss \ decomposition
\begin{equation}
P\left(  x\right)  =h_{0}\left(  x\right)  +\left|  x\right|  ^{2}h_{1}\left(
x\right)  +...+\left|  x\right|  ^{2N}h_{N}\left(  x\right)  . \label{Gauss23}%
\end{equation}
Then
\begin{equation}
d\left(  P\cdot Y_{k}\right)  \leq\max_{r=0,...,N}\left\{  r+\deg
h_{r}\right\}  \leq\deg P\left(  x\right)  . \label{sharp}%
\end{equation}

\end{corollary}

%

\proof
By (\ref{additiv}) $d\left(  \left|  x\right|  ^{2r}h_{r}Y_{k}\right)
=r+d\left(  h_{r}Y_{k}\right)  .$ By Proposition \ref{P23} $d\left(
h_{r}Y_{k}\right)  \leq\min\left\{  \deg h_{r},\deg Y_{k}\right\}  \leq\deg
h_{r}.$ This proves the first inequality. Further we know that $\deg\left|
x\right|  ^{2r}h_{r}=2r+\deg h_{r}\leq\deg P$ for $r=0,...,N$. Hence the
second inequality is established.%
\endproof

In the following we want to give an explicit formula for $N_{P}.$

\begin{theorem}
\label{T39}Let $Y_{k,m}\left(  x\right)  $ be an orthonormal basis of
spherical harmonics with $k\in\mathbb{N}_{0}$ and $m=1,...,a_{k}.$ Then
$d\left(  Y_{k,m}\left(  x\right)  Y_{k,m_{1}}\left(  x\right)  \right)  =k$
if and only if $m=m_{1}.$
\end{theorem}

%

\proof
We start with a general remark: \label{minkl}Let $Y_{k}$ and $Y_{l}$ be
harmonic homogeneous polynomials of degree $k$ and $l$ respectively. Clearly
$Y_{k}\left(  x\right)  Y_{l}\left(  x\right)  $ is a homogeneous polynomial
of degree $k+l.$ By Proposition \ref{P23} it has polyharmonic degree at most
$\min\left\{  k,l\right\}  .$ By Gau\ss \ decomposition there exist harmonic
homogeneous polynomials $h_{k+l-2u},$ either $h_{k+l-2u}$ is zero or of exact
degree $k+l-2u$ for $u=0,...,\min\left\{  k,l\right\}  ,$ such that
\begin{equation}
Y_{k}\left(  x\right)  Y_{l}\left(  x\right)  =\sum_{u=0}^{\min\left\{
k,l\right\}  }\left|  x\right|  ^{2u}h_{k+l-2u}(x). \label{last1}%
\end{equation}
Now assume that $Y_{k}\left(  x\right)  =Y_{k,m}\left(  x\right)  $ and
$Y_{l}\left(  x\right)  =Y_{k,m_{1}}\left(  x\right)  .$ Let us consider the
summand $\left|  x\right|  ^{2k}h_{0}\left(  x\right)  $ for $u=k.$ Then
$h_{0}$ must have degree $0,$ hence it is a constant polynomial. Integrate
equation (\ref{last1}) with respect to $d\theta.$ Since $h_{2k-2u}$ is either
$0$ or of exact degree $2k-2u>0$ for $u=1,...,k$ the integral over the sphere
of $\left|  x\right|  ^{2u}h_{k+l-2u}(x)$ will vanish. Then we obtain
\[
\delta_{m,m_{1}}\left|  x\right|  ^{2k}=\int_{\mathbb{S}^{n-1}}h_{0}%
d\theta=h_{0}\omega_{n}.
\]
Hence for $m\neq m_{1}$ we see that the polyharmonic degree is less than $k,$
for $m=m_{1}$ it is exactly $k.$ The proof is finished.%
\endproof

\begin{theorem}
\label{T40}Let $P\left(  x\right)  $ be a homogeneous polynomial of degree
$N,$ say of the form
\[
P\left(  x\right)  =\sum_{t,k\in\mathbb{N}_{0},2t+k=N}\sum_{m=1}^{a_{k}%
}a_{t,k,m}\left|  x\right|  ^{2t}Y_{k,m}\left(  x\right)  .
\]
Let $k_{0}$ be the largest natural number such that $a_{t_{0},k_{0},m_{0}}%
\neq0$ for some $m_{0}$ in the above sum. Then
\[
N_{P}:=\sup_{k\in\mathbb{N}_{0},m=1,...,a_{k}}d\left(  P\left(  x\right)
Y_{k,m}\left(  x\right)  \right)  =\frac{1}{2}\left(  N+k_{0}\right)  .
\]

\end{theorem}

%

\proof
Since $d\left(  P+Q\right)  \leq\max\left\{  d\left(  P\right)  ,d\left(
Q\right)  \right\}  $ we obtain for $k_{1}\in N_{0}$ and $m_{1}\in\left\{
1,...,a_{k_{1}}\right\}  $ that
\[
d\left(  P\left(  x\right)  Y_{k_{1},m_{1}}\left(  x\right)  \right)  \leq\max
d\left(  \left|  x\right|  ^{2t}Y_{k,m}Y_{k_{1},m_{1}}\left(  x\right)
\right)
\]
where the maximum ranges over all indices $t,k,m$ with $a_{t,k,m}\neq0.$ Since
$d\left(  Y_{k,m}Y_{k_{1},m_{1}}\right)  \leq k$ we arrive at (note that
$2t+k=N)$
\[
d\left(  P\left(  x\right)  Y_{k_{1},m_{1}}\left(  x\right)  \right)  \leq
\max\left\{  t+k\right\}  =\frac{1}{2}\max\left\{  N+k\right\}  \leq\frac
{1}{2}\left(  N+k_{0}\right)  .
\]
Hence we see that $\frac{1}{2}\left(  N+k_{0}\right)  $ is a bound for the
polyharmonic degree of $P\left(  x\right)  Y_{k_{1},m_{1}}\left(  x\right)  .$

Let us consider $P\left(  x\right)  Y_{k_{0},m_{0}}\left(  x\right)  $ where
$k_{0}$ is as in the theorem. Consider a summand $a_{t,k,m}\left|  x\right|
^{2t}Y_{k,m}$ with $a_{t,k,m}\neq0.$ Then $k\leq k_{0}$ and Proposition
\ref{P23} shows that $d\left(  Y_{k,m}Y_{k_{0},m_{0}}\right)  \leq k,$ hence
for $k<k_{0}$ each summand $a_{t,k,m}\left|  x\right|  ^{2t}Y_{k,m}%
Y_{k_{0},m_{0}}$ has polyharmonic degree
\begin{equation}
d\left(  a_{t,k,m}\left|  x\right|  ^{2t}Y_{k,m}Y_{k_{0},m_{0}}\right)  \leq
t+k=\frac{1}{2}\left(  N+k\right)  <\frac{1}{2}\left(  N+k_{0}\right)  .
\label{holds}%
\end{equation}
Now consider the case $k=k_{0}.$ If $m\neq m_{0}$ then we apply Theorem
\ref{T39} and the same argument shows that (\ref{holds}) holds. Finally assume
that $k=k_{0}$ and $m=m_{0}.$ Then Theorem \ref{T39} shows that $a_{t_{0}%
,k_{0},m_{0}}\left|  x\right|  ^{2t_{0}}Y_{k_{0},m_{0}}Y_{k_{0},m_{0}}$ has
exact polyharmonic degree $t_{0}+k_{0}=\frac{1}{2}\left(  N+k_{0}\right)  . $
Hence we have proven that
\[
P\left(  x\right)  Y_{k_{0},m_{0}}=a_{t_{0},k_{0},m_{0}}\omega_{n}\left|
x\right|  ^{N+k_{0}}+R\left(  x\right)
\]
where $R\left(  x\right)  $ has polyharmonic degree $<\frac{1}{2}\left(
N+k_{0}\right)  $. Thus $P\left(  x\right)  Y_{k_{0},m_{0}}$ has exact
polyharmonic degree $\frac{1}{2}\left(  N+k_{0}\right)  .$
\endproof

Let us finish with the following remark. Let $P\left(  x\right)  $ be an
arbitrary polynomial. We can write $P\left(  x\right)  =\sum_{j=0}^{N}%
P_{j}\left(  x\right)  $ where $P_{j}\left(  x\right)  $ are homogeneous
polynomials. It is not very difficult to see that
\[
d\left(  P\cdot Y_{k,m}\right)  =\max_{j=0,...,N}d\left(  P_{j}\cdot
Y_{k,m}\right)  ,
\]
see e.g. the proof of Theorem 1.27 in \cite{ABR92}. Hence $N_{P}$ is the
maximum of $N_{P_{j}}$ for $j=0,...,N.$

1. Ognyan Kounchev, Institute of Mathematics and Informatics, Bulgarian
Academy of Sciences, 8 Acad. G. Bonchev Str., 1113 Sofia, Bulgaria;

e--mail: kounchev@math.bas.bg, kounchev@math.uni--duisburg.de

2. Hermann Render, Departamento de Matem\'{a}ticas y Computati\'{o}n,
Universidad de la Rioja, Edificio Vives, Luis de Ulloa, s/n. 26004
Logro\~{n}o, Spain; e-mail: render@gmx.de; hermann.render@dmc.unirioja.es

\end{document}